\documentclass{article}
\usepackage{comment}
\usepackage{amsfonts}
\usepackage{bm}
\usepackage{diagbox}
\usepackage{color}
\usepackage{amsmath}
\usepackage{amsfonts}
\usepackage{amsthm}
\usepackage{graphicx}

\theoremstyle{plain} 
\newtheorem{theorem}{Theorem}[section]

\theoremstyle{definition} 

\theoremstyle{remark} 
\newtheorem{remark}[theorem]{Remark}

\title{A BFBt preconditioner for Double Saddle-Point Systems \\ 
{\it \large In memory of Howard Elman}}

\author{%
	{\sc
		Chen Greif\thanks{Email: greif@cs.ubc.ca. The work of the author was supported in part by the Natural Sciences and Engineering Research Council of Canada.}} \\[2pt]
		Department of Computer Science, The University of British Columbia\\
		Vancouver, Canada
}
\date{\em To be published in  IMA Journal of Numerical Analysis}

\begin{document}

\maketitle

\begin{abstract}
{We consider block preconditioners for double saddle-point systems, and investigate the effect of approximating the nested Schur complement associated with the trailing diagonal block on the eigenvalue distribution of the preconditioned matrix. We develop a variant of Elman's BFBt method
and adapt it to this family of linear systems. Our findings are illustrated on a Marker-and-Cell discretization of the Stokes-Darcy equations.}
{double saddle-point system; preconditioning; BFBt; least-squares commutator; Schur complement; Stokes-Darcy equations}
\end{abstract}

\section{Introduction} \label{sec:intro}

Consider large and sparse double saddle-point systems of the form
\begin{equation}
	\underbrace{\begin{bmatrix}
			A&{{B^T}}&0\\
			B&-D&{{C^T}}\\
			0&C&0
	\end{bmatrix}}_{\mathcal K} 
\begin{bmatrix}
			u\\
			v\\
			w
	\end{bmatrix}=
\begin{bmatrix}
			b\\
			c\\
			d
	\end{bmatrix},
\label{eq:system}
\end{equation}
where $A\in \mathbb{R}^{n\times n}, B\in \mathbb{R}^{m\times n}$,  $C\in \mathbb{R}^{p\times m}$, and $D\in \mathbb{R}^{m\times m}$.  
We assume that   $\mathcal{K}$ is nonsingular; see \cite[]{xie2020note, beik2024invertibility} 
 for necessary and sufficient conditions for invertibility. 
Our interest  is in minimal residual iterative solvers that use Schur complements for preconditioning, and for that purpose we will assume that the following three matrices are well defined and nonsingular: 
\begin{itemize}
\item The leading block, $A$;
\item the (negated) Schur complement, \begin{equation} 
S_1=D+B A^{-1} B^T  ;
\label{eq:S1}
\end{equation}
\item the nested Schur complement,
\begin{equation}
S_2=C S_1^{-1} C^T.
\label{eq:S2}
\end{equation}
\end{itemize}
The algebraic relationship between $n$ and $m$ may vary, and it depends on the specific application. 
Notice that if ${\rm rank}(C)=p \le m$ and $S_1$ exists and is nonsingular, then so is $S_2$.  

If $A$ and/or $D$ are nonsymmetric and $\mathcal{K}$ is not too far from being normal, a block LDU decomposition provides a suitable potential candidate for a block-triangular preconditioner: we will study the block lower-triangular matrix
\begin{equation}
	\mathcal{M}_{LT}=\begin{bmatrix}
			A & 0 & 0 \\
			B & -S_1 & 0 \\
			0 & C  & S_2
	\end{bmatrix}.
\label{eq:pc1}
\end{equation}

An alternative, structurally symmetric preconditioning approach, is the block-diagonal preconditioner
\begin{equation}
	\mathcal{M}_D=\begin{bmatrix}
			A & 0 & 0 \\
			0 & S_1 & 0 \\
			0 & 0  & S_2
	\end{bmatrix}.
\label{eq:pc2}
\end{equation}

In practice, $\mathcal{M}_{LT}$ and $\mathcal{M}_D$ are not applied as preconditioners in the forms  \eqref{eq:pc1} and \eqref{eq:pc2}, respectively, because  the Schur complements $S_1$ and $S_2$ are typically prohibitively expensive to compute or invert. Nonetheless,  observations on the structure of the preconditioned eigenvalues are useful because they serve as a basis for seeking preconditioners that rely on using effective approximations of  $S_1$ and $S_2$ or the action of their inverses;  this is similar to the state of affairs with (classical) saddle-point systems \cite[]{bgl05, ESW2014}.

Several recent papers consider preconditioners of the forms \eqref{eq:pc1} and \eqref{eq:pc2} and other variants, where $A, S_1$, and $S_2$ are replaced by approximations
$\widehat{A}$, $\widehat{S}_1$ and $\widehat{S}_2$, respectively. Those papers provide theoretical and experimental insights on the spectra of the preconditioned matrices and the convergence of corresponding iterative solvers. In \cite[]{hm19}, the authors consider $\widehat{A} \approx A$, $\widehat{S}_1 \approx S_1$, and $\widehat{S}_2 = C \widehat{S}_1^{-1} C^T \approx S_2$.  Several other papers, e.g., \cite[]{Balani, Balani2}, go further in the direction of developing effective approximations, and demonstrate theoretically and experimentally the viablity of the approaches taken. Eigenvalue bounds for preconditioned symmetric double saddle-point systems are derived in \cite{BG2023}. In \cite{sogn18}, important observations on the eigenvalues of symmetric multiple saddle-point systems are given, and in \cite[]{bergamaschi2024spectral, bergamaschi2025eigenvalueboundspreconditionedsymmetric, pearson2023double, pp2023},  scalable preconditioners for such systems are developed and eigenvalue bounds are derived. 
 Extensions to block preconditioners that do not have a block-diagonal form have also been analyzed, for example, in \cite[]{Balani, xie2020note}.  Modified double-saddle point systems---for example, where one of $B$ and $B^T$ is negated and similarly for $C$, giving rise to a nonsymmetric block matrix whose symmetric part is positive semidefinite---have also been considered  \cite[]{hm19}, using shift-splitting preconditioners \cite[]{cao2019} and other approaches \cite[]{Balani}.

A question we are interested in is what happens when $A$ and $S_1$ are solved for exactly and $S_2$ is given approximately.  We show that in this case the eigenvalues of $\widehat{S}_2^{-1} S_2$ play a central role in the convergence of minimal residual preconditioned iterative solvers. While this question may seem to have a limited theoretical value, it is potentially useful, because it shows that some theoretical observations on the convergence of iterative solvers  for the indefinite $(n+m+p) \times (n+m+p)$ linear system \eqref{eq:system} can be inferred by studying a smaller generalized eigenvalue problem associated with $p \times p$ matrices.

In the context of finding effective approximations for $S_2$, a primary goal in this paper is to develop a variant of Elman's BFBt method for double saddle-point systems. This method has led to the popular least-squares commutator  \cite[]{Elman99, ehsst2006, ehsst2008}, and was designed for nonsymmetric saddle-point systems arising from the Navier-Stokes equations. It amounts to a block-triangular preconditioner. We adapt it to double saddle-point systems, and demonstrate its effectiveness on the Marker-and-Cell discretization of the Stokes-Darcy equations, using practical approximations of $S_1$ and~$S_2$. 

In Section \ref{sec:inexS2} we perform some basic eigenvalue analysis, exploring the effect that using an approximation of $S_2$ has on the block preconditioners. In Section \ref{sec:LSC} we consider a BFBt approximation of $S_2^{-1}$ for double saddle-point systems and provide some additional spectral analysis. In Section \ref{sec:numerical} we explore the  Stokes-Darcy equations as a test case and present some numerical results.  Finally, in Section \ref{sec:conc} we offer brief concluding remarks.

\section{An inexact $S_2$}
\label{sec:inexS2}

Suppose we compute exact solutions of linear systems involving $A$ and $S_1$, while $S_2$ is replaced by a  matrix $\widehat{S}_2$ that is easier to invert.
We denote the corresponding block lower-triangular preconditioner based on 
$\mathcal{M}_{LT}$
 by
\begin{equation}
\mathcal{\widehat{M}}_{LT}  = \begin{bmatrix}
A & 0 &0\\
			B & -S_1 & 0 \\
			0&C &\widehat{S}_2
	\end{bmatrix},
\label{eq:hatMLT}
\end{equation}
and the block-diagonal preconditioner based on $\mathcal{M}_D$ 
by
\begin{equation}
\mathcal{\widehat{M}}_{D}  = \begin{bmatrix}
A & 0 &0\\
			0 & S_1 & 0 \\
			0&0 &\widehat{S}_2
	\end{bmatrix}.
\label{eq:hatMLTD}
\end{equation}

\subsection{Eigenvalues of $\mathcal{\widehat{M}}_{LT}^{-1} \mathcal{K}$}
The preconditioned matrix is given by
\begin{equation}
\mathcal{\widehat{M}}_{LT}^{-1} \mathcal{K}  = \begin{bmatrix}
I_n & A^{-1} B^T &0\\
			0 & I_m & -S_1^{-1} B^T \\
			0&0 &\widehat{S}_2^{-1} S_2
	\end{bmatrix}.
\label{eq:hatMLTK}
\end{equation}
 When $\widehat{S}_2=S_2$,  the eigenvalues of the preconditioned matrix $\mathcal{M}_{LT}^{-1} \mathcal{K}$ are all equal to 1, with a minimal polynomial of degree 3 \cite[Theorem 4.4]{GH2023}.
If $\widehat{S}_2 \ne S_2$ it readily follows that $n+m$ eigenvalues equal to 1 are revealed from the first two leading blocks of the preconditioned matrix, and the distribution of the remaining eigenvalues of the preconditioned matrix depends on the generalized eigenvalues
\begin{equation}
  S_2 z = \mu  \widehat{S}_2 z. 
\label{eq:geneig}
\end{equation}

\subsection{Eigenvalues of $\mathcal{\widehat{M}}_{D}^{-1} \mathcal{K}$} 
\label{sec:symcase}

The eigenvalue problem corresponding to the eigenvalues of the preconditioned matrix $\mathcal{\widehat{M}}_{D}^{-1} \mathcal{K}$ is
\begin{equation}
\underbrace{
\begin{bmatrix}
A^{-1} & 0 &0\\
			0 & S_1^{-1} & 0 \\
			0&0 &\widehat{S}_2^{-1}
	\end{bmatrix}}_{{\mathcal{\widehat{M}}}_D^{-1}}
	\underbrace{
\begin{bmatrix}
			A&{{B^T}}&0\\
			B& -D &{{C^T}}\\
			0&C&0
	\end{bmatrix}}_{\mathcal{K}}
\begin{bmatrix}
			x\\
			y\\
			z
	\end{bmatrix}=\lambda
\begin{bmatrix}
			x\\
			y\\
			z
	\end{bmatrix},
\label{eq:eigMK0}
\end{equation}
which gives us the  matrix equations
\begin{subequations}
\begin{align}
x + A^{-1} B^T y = \lambda x;  \label{eq:A0} \\
S_1^{-1} B x  - S_1^{-1} D y + S_1^{-1} C^T z = \lambda y ; \label{eq:B0} \\
 \widehat{S}_2^{-1} C y = \lambda z.  \label{eq:C0}
 \end{align}
\label{eq:ABC}
\end{subequations}
We discuss two common cases.

\subsubsection{$D \ne 0$ and $m \ge n$} 
\label{sec:Dne0}
The block-tridiagonal structure that we assume in \eqref{eq:system} is not necessarily the only way to define a double saddle-point system. In \cite{beik18}, for example, the matrix appears in a form where the (1,3) block is nonzero and the (2,3) block is zero. A symmetric block permutation can be applied to transform the system into the block-tridiagonal form \eqref{eq:system}. This flips the algebraic relationship between the dimensions of the (1,1) block and  the (2,2) block. 

The Marker-and-Cell discretization of the Stokes-Darcy equations is a relevant instance. Here, the matrix $D$ is larger in dimensions than $A$  and $B$ is rank deficient.
An eigenvalue analysis for $\widehat{S}_2=S_2$ in \eqref{eq:eigMK0} and \eqref{eq:C0} was performed in \cite[Theorem 4.1]{GH2023}, and it was
shown that $\lambda=1$  and $\lambda=-1$ are eigenvalues with algebraic multiplicities equal or nearly equal $n$, along with other observations.

Consider now the  case of interest to us,  $\widehat{S}_2 \ne S_2$. Suppose $x=0$ but $y \ne 0$ and $z \ne 0$. 
From \eqref{eq:A0} we have 
$ B A^{-1} B^T y = (S_1-D) y= 0,$
and from \eqref{eq:B0}, multiplied by $S_1$, we have
$  -Dy + C^T z = -S_1 y + C^T z = \lambda S_1 y,$ which yields
$  y = \frac{1}{\lambda+1} S_1^{-1} C^T z.$
Substituting this in \eqref{eq:C0}, we get
$$ \widehat{S}_2^{-1} S_2 z = \lambda (\lambda+1) z.$$
Denoting, as in \eqref{eq:geneig}, the eigenvalues of $ \widehat{S}_2^{-1} S_2$ by $\mu$, we have the relation
$$ \lambda(\lambda+1) = \mu.$$
Therefore
$$ \lambda_{\pm} = \frac{-1 \pm \sqrt{1+4 \mu}}{2}.$$
If $\mu=1$, then we get $\lambda_{\pm} = \frac{-1 \pm \sqrt{5}}{2}.$ The algebraic multiplicities of this pair of eigenvalues is equal to the multiplicity of $\mu=1$ for the generalized eigenvalue problem \eqref{eq:geneig}.

\subsubsection{$A$ is symmetric positive definite, $D = 0$, and $n \ge m$} 
\label{sec:symD0}
This case has been extensively analyzed in the literature. The nonsingularity of $\mathcal{K}$ implies that ${\rm rank}(B)=m$. 
The Schur complements $S_1$ and $S_2$ are symmetric positive definite, and $\mathcal{K}$ is symmetric indefinite.

The preconditioned matrix $\mathcal{M}_D^{-1} \mathcal{K}$ is known to have precisely six distinct eigenvalues \cite[]{sogn18}, which can be compactly written as 
\begin{equation*}
\lambda_{i,j} = 2 \cos \left( \frac{2i+1}{2j+3} \pi \right), \qquad  j=0, 1, 2, \quad i=0, 1,\dots, j.
\end{equation*}
Explicitly, these six eigenvalues values are given (to the precision displayed) by $$\{ -1.2470, -0.6180, 0.4450, 1, 1.6180, 1.8019 \}. $$
See also \cite[]{BG2023, cai2021schur,  pearson2023double}.
In the absence of roundoff errors and if we assume we have exact solvers for linear systems involving $A$, $S_1$ and $S_2$, a minimal residual iterative scheme such as MINRES \cite[]{ps1975}, which uses $\mathcal{M}_D$ as a preconditioner for solving the linear system \eqref{eq:system}, is expected to converge within six iterations. 

\begin{theorem}Suppose $A$ is symmetric positive definite, $D=0$, and $n \ge m$, and suppose further that $\widehat{S}_2$ is symmetric positive definite. 
Then, the eigenvalues of the preconditioned matrix ${{\mathcal{\widehat{M}}}_{D}^{-1}} \mathcal{K}$ defined in \eqref{eq:eigMK0} are:
\begin{itemize}
\item $\lambda=1$, of algebraic multiplicity $n-m$, with corresponding eigenvector $(x,0,0)$, where $x \in {\rm ker}(B)$ is a set of $n-m$ linearly-independent null vectors of $B$ (a basis for the kernel of $B$).
\item $\lambda=\frac{1 \pm \sqrt{5}}{2}$ ($\lambda_+ \approx 1.618$ and $\lambda_- \approx -0.618$), each of algebraic multiplicity $m-p$, with corresponding eigenvectors $(x,y,0)$ that satisfy $0 \ne y \in {\rm ker}(C)$ and $x=\frac{1}{\lambda-1} A^{-1} B^T$. 
\item
The remaining $3p$ eigenvalues are the roots of the $p$ cubic polynomials 
$$\lambda^3 - \lambda^2 - (1+\mu_i) \lambda+\mu_i=0,  $$
where $\mu_i, \ i=1,\dots,p$, are the $p$ eigenvalues of the generalized eigenvalue problem~\eqref{eq:geneig}.

\end{itemize}
\label{thm:eig}
\end{theorem}

\begin{proof}
We consider \eqref{eq:ABC} and split the discussion into a few separate cases.

\noindent
{\em The eigenvalue $\lambda=1$.} Multiplying \eqref{eq:A0} by $B$, we  observe that for any null vector $x \in {\rm ker}(B)$ $\lambda=1$ is an eigenvalue with eigenvectors $(x,0,0)$. It is immediate to confirm that there is no situation in which $\lambda=1$ is an eigenvalue with a corresponding eigenvector $(x,y,z)$ for which $x \in {\rm ker}(B)$ and $y \ne 0$ or $z \ne 0$. Suppose now that $x \notin {\rm ker}(B)$. From \eqref{eq:A0}, multiplied by $B$, we have that $S_1 y=(\lambda-1) Bx$ with the assumption  $Bx\ne 0$. But if $\lambda=1$ then we still have that $S_1y=0$, which implies that $y=0$, and this also leads to $z=0$. However, \eqref{eq:B0} in this case would imply $Bx=0$, which is a contradiction. Therefore,   we have that for $\lambda=1$ we must have
$x \in {\rm ker}(B)$, and then necessarily $y=z=0$ and $(x,0,0)$ is an eigenvector. We have thus identified all possible cases for $\lambda=1$, and we can conclude that its algebraic multiplicity is precisely $n-m$, the nullity of $B$. For $x$ we can take any set of $n-m$ linearly-independent null vectors of $B$. 

\noindent
{\em The eigenvalues $\lambda_{\pm}= \frac{1 \pm \sqrt{5}}{2}$.}
Assume that $x \notin {\rm ker}(B)$ and $\lambda \ne 1.$ Starting with \eqref{eq:C0}, suppose $y \in {\rm ker}(C)$. Notice that we cannot have  $y=0$ because it immediately leads to $x=z=0$. Given our rank assumptions, namely that $C$ has full row rank, there are $m-p$ linearly-independent null vectors $y$ of $C$. From \eqref{eq:C0}, since $\lambda \ne 0$ we have $z=0$. This allows us to combine \eqref{eq:A0}--\eqref{eq:B0} and obtain $$ \lambda^2-\lambda-1=0,$$ from which we get the golden ratio eigenvalues $\lambda_{\pm}= \frac{1 \pm \sqrt{5}}{2},$ each associated to  $m-p$  eigenvectors
$(x,y,0)$ that satisfy $x \notin {\rm ker}(B)$ and $0 \ne y \in {\rm ker}(C)$.  Once $y$ has been determined (as a null vector of $C$), the vector $x$ is determined by \eqref{eq:A0}
as $x=\frac{1}{\lambda-1} A^{-1} B^T y.$ We will show below that it is not possible to have additional instances of the same eigenvalues with $y \notin {\rm ker}(C)$, and therefore the algebraic multiplicity of each of these two eigenvalues is precisely $m-p$. If $m=p$, there are no eigenvalues in this category. \\

\noindent
{\em Remaining eigenvalues.} In total, we have accounted so far for $n-m+2(m-p)=n+m-2p$ eigenvalues and their corresponding eigenvectors. There are $3p$ more eigenpairs to determine. We now assume that the remaining eigenvectors satisfy $x \notin {\rm ker}(B)$ and $ y \notin {\rm ker}(C)$.
Multiplying \eqref{eq:A0} by $B$ and \eqref{eq:B0} by $C$ and combining the two resulting equations, we obtain
$$ \left(\lambda - \frac{1}{\lambda-1} \right) C y = S_2 z.$$
Note that we have already established that for this case $\lambda \ne 1$, and therefore there is no concern of division by zero. 
Substituting  $ Cy = \lambda \widehat{S}_2 z $ from \eqref{eq:C0},
we obtain  the generalized eigenvalue problem
\begin{equation} S_2 z = \underbrace{\lambda \left(\frac{\lambda^2-\lambda-1}{\lambda-1} \right)}_\mu  \widehat{S}_2 z.
\label{eq:mu}
\end{equation}
Since both $\widehat{S}_2$ and $S_2$ are assumed to be symmetric positive definite, all the generalized eigenvalues $\mu$ are positive, and therefore we cannot have $\lambda^2-\lambda-1=0$, which shows that we cannot have any additional golden ratio eigenvalues. Equating $\mu$ with the expression for $\lambda$ on the right-hand side of \eqref{eq:mu} simplifies to the cubic equation for $\lambda$ given in the statement of the theorem, as required.
\end{proof}

\begin{remark}
{\it If $\widehat{S}_2=S_2$ in \eqref{eq:mu}, then 
$ \lambda^3-\lambda^2-2 \lambda+1=0$, $\mathcal{\widehat{M}}_D \equiv \mathcal{{M}}_D$,  and the corresponding eigenvalues are given by 
$ \lambda_1= 2 \cos \left(\frac{ \pi}{7} \right) \approx  1.8019; \ \lambda_2 =  2 \cos \left(\frac{3 \pi}{7} \right)  \approx  0.4450; \ \lambda_3 = 
2 \cos \left(\frac{5 \pi}{7} \right) \approx  -1.2470,$
each of multiplicity $p$,
which leads to the six distinct eigenvalues of ${\mathcal{M}}_D^{-1} \mathcal{K}$  mentioned earlier in this section.}
\label{rem:cubic}
\end{remark}

\subsubsection{Perturbation analysis}
In the case described in  Section \ref{sec:symD0}, the generalized eigenvalues of \eqref{eq:mu} satisfy
\begin{equation}
 \lambda^3-\lambda^2 -(1+\mu)\lambda + \mu=0.
\label{eq:mulam}
\end{equation}
In Section \ref{sec:LSC} we will consider in detail a situation where the generalized eigevalues satisfy $\mu \ge 1$.  In that case, we have
\begin{equation}
1 \leq \mu = \lambda^2-\frac{\lambda}{\lambda-1}.
\label{eq:mu2}
\end{equation}
Denoting $\lambda := \lambda(\mu)$, the roots stated in Remark~\ref{rem:cubic} are $\lambda_1(1), \lambda_2(1)$, and $\lambda_3(1)$.
Notice that $\lambda \ne 1$, consistently with \eqref{eq:mulam}. 

Following \cite[Lemma 2.2]{bergamaschi2024spectral}, we differentiate $\lambda=\lambda(\mu)$ with respect to $\mu$ and obtain  
$$\lambda'=\left(2\lambda + \frac{1}{(\lambda-1)^2} \right)^{-1}.$$%
From this it follows that $\lambda'>0$ for $\lambda>0$, which means that $\lambda_1(\mu)$ and $\lambda_2(\mu)$ are increasing for small perturbations of $\mu$ deviating from 1. On the other hand, for $\lambda<-1$, we have that $\lambda'<0$, and therefore $\lambda_3(\mu)$ is decreasing. From this it follows that the perturbed eigenvalues (under sufficiently small perturbations)  move away from the origin and are in the intervals 
$$[ \lambda_3(\mu_{\max}), -1.2470], \qquad [0.4450, 1), \qquad [1.8019,\lambda_1(\mu_{\max})].$$

In order to assess the MINRES convergence rate in this case, it is necessary to express $\lambda_1$ and $\lambda_3$ more explicitly in terms of $\mu_{\max}$. 
Suppose $\mu=1+\varepsilon$, where $\varepsilon \ll 1.$ Then we can write the asymptotic expansion
$$ \lambda = \lambda^{(0)} + \varepsilon \lambda^{(1)} + \mathcal{O}(\varepsilon^2)$$ 
and substitute it in \eqref{eq:mulam}.
The values for $\lambda^{(0)}$ are the original roots stated in Remark~\ref{rem:cubic}, and 
$$ \lambda^{(1)} = \frac{\lambda^{(0)}-1}{3 (\lambda^{(0)})^2 - 2 \lambda^{(0)} -2}.$$
A short calculation shows that for $\lambda^{(0)}=1.8019$ we have  $\lambda^{(1)}=0.1938$ and for $\lambda^{(0)}=-1.2470$ we have $\lambda^{(1)}=-0.4356$. In either case the $\mathcal{O}(\varepsilon)$ term in the asymptotic expansion is moderate in magnitude, which tells us that up to first-order, small perturbations in $\mu$ result in perturbations of a similar magnitude in $\lambda$, and the intervals for MINRES can be characterized accordingly.

In practice, $\mu_{\max}$ may not be close to 1. In this case, the perturbation analysis is less useful, but we may still be able to establish that the preconditioner is  effective, especially if many of the generalized eigenvalues are close to 1, since the eigenvalue intervals specified above are bounded uniformly away from zero.  
If $\mu_{\max} \gg 1$, then since $\lambda^2$ dominates under this scenario  the next term in \eqref{eq:mu2}, we have 
$$ \lambda_1(\mu_{\max}) \approx  \sqrt{\mu_{\max}}, \qquad  \lambda_3(\mu_{\max}) \approx -\sqrt{\mu_{\max}}.$$

\section{BFBt preconditioning}
\label{sec:LSC}

In \cite[]{Elman99}, Elman considered  saddle-point matrices  arising from the numerical solution of Navier-Stokes equations
$$\mathcal{K} = \begin{bmatrix} F & G^T \\ G & 0 \end{bmatrix},$$
where $F$ is a (nonsymmetric) discrete convection-diffusion operator and $G$ is a discrete divergence operator. Preconditioners of the form
$$\mathcal{M} = \begin{bmatrix} F & G^T \\ 0 & -X\end{bmatrix}$$
were proposed,
where $X$ is an approximation of the inverse of the  Schur complement, $S=G F^{-1} G^T$.
Spectral and range-space considerations were used to derive the approximation 
$$ X = (G G^T)  (G F G^T)^{-1} (G G^T).$$
A key feature here is that upon inversion, in $X^{-1}$ the matrix $F$ is not inverted and the operation of the inverse of $G G^T$ (or scaled variants) is relatively inexpensive for the Navier-Stokes equations, since it is related to the scalar Laplacian. 

This preconditioner was termed {\em BFBt}  in \cite[]{Elman99}.  
It was later further tailored to finite element discretizations of the Navier-Stokes equations by using a scaling matrix \cite[]{ehsst2006, ehsst2008}: replace $G G^T$ by $G \widehat{Q}^{-1} G^T$, where $\widehat{Q}$ is the diagonal of the scalar mass matrix. The resulting preconditioner was termed the {\em least-squares commutator} and has been shown to be a highly effective preconditioning approach for the finite element discretization of the Navier-Stokes equations; see, e.g., \cite[]{ESW2014}.

\subsection{Application to double saddle-point systems}
 For the double saddle-point system \eqref{eq:system}, suppose inverting $C C^T$ is relatively cheap, and  consider the  choice
\begin{equation}
\widehat{S}_2^{-1} = (C C^T)^{-1} C S_1 C^T (C C^T)^{-1}. 
\label{eq:bfbt}
\end{equation}

 In Section \ref{sec:inexS2} we showed that some of the eigenvalues of the preconditioned matrices we consider  depend directly on the generalized eigenvalue problem \eqref{eq:geneig}.
For the choice \eqref{eq:bfbt}, we have
\begin{equation}
 \widehat{S}_2^{-1} S_2 =(C C^T)^{-1} C S_1 \underbrace{C^T (C C^T)^{-1} C}_{P} S_1^{-1} C^T.
\label{eq:S2m1S2}
\end{equation}
The orthogonal projector $P=P^2$ behaves like the identity on vectors in the range space of $C^T$.

\subsection{Eigenvalue analysis}

Suppose $m > p$; the special case $m=p$  is omitted.
We perform two main steps, applying and extending an elegant mathematical trick  used in \cite[]{Elman99}:
\begin{enumerate} \item {\bf Step 1.} We express the spectrum of $\widehat{S}_2^{-1} S_2$ in terms of $S_1$ and the singular values \& vectors of the ``tall and skinny'' matrix $C^T$, and obtain via a similarity transformation a matrix that is slightly simpler than $\widehat{S}_2^{-1} S_2$.
\item 
{\bf Step 2.} We express the simpler matrix by partitioning a related matrix into a $2 \times 2$ block matrix, finding its inverse, and extracting its leading block. This  turns out to be useful in  identifying algebraic multiplicities of some of the eigenvalues of the preconditioned matrix. 
\end{enumerate}

\subsubsection{A similarity transformation}

We start with Step 1, as follows. Consider the singular value decomposition (SVD) of $C^T$:
$$  C^T = U \Sigma V^T, $$
where $U$ is $m \times m$ and orthogonal, $\Sigma$ is $m\times p$ and  contains the singular values of $C^T$ along the diagonal of its upper $p \times p$ part, and $V$ is $p \times p$ and orthogonal.
Let us establish some additional notation, as follows: $$ \Sigma = \begin{bmatrix} \widehat{\Sigma} \\ 0 \end{bmatrix} \qquad {\rm and} \qquad 
U = \begin{bmatrix}\widehat{U} & \widetilde{U} \end{bmatrix}, $$
where $\widehat{\Sigma}$ is a square diagonal $p \times p$ matrix, $\widehat{U}$ is $m \times p$ with orthonormal columns, and $\widetilde{U}$ is $m \times (m-p)$ with orthonormal columns. Accordingly, $$C^T = \widehat{U} \widehat{\Sigma} V^T$$ is the reduced SVD of $C^T$.

Writing down $\widehat{S}_2$ and $S_2$ in terms of $S_1$ and the reduced singular vectors and singular values of $C^T$, and  substituting the SVD into \eqref{eq:S2m1S2}, we obtain
$$ \widehat{S}_2^{-1} S_2 = V \widehat{\Sigma}^{-1} (\widehat{U}^T S_1 \widehat{U}) (\widehat{U}^T S_1^{-1} \widehat{U}) \widehat{\Sigma} V^T. $$
From this it follows that $\widehat{S}_2^{-1} S_2$ is similar to 
\begin{equation}
W = (\widehat{U}^T S_1 \widehat{U}) (\widehat{U}^T S_1^{-1} \widehat{U}), 
\label{eq:W}
\end{equation}
and therefore these two matrices have the same eigenvalues. 

We now proceed to Step 2. We observe that the two matrices whose product forms the matrix $W$ can be extracted from a $2 \times 2$ partitioning, as follows. On the one hand, 
\begin{equation} (U^T S_1 U)^{-1} = \begin{bmatrix} \widehat{U}^T S_1 \widehat{U} & \widehat{U}^T S_1 \widetilde{U} \\
\widetilde{U}^T S_1 \widehat{U} \ & \widetilde{U}^T S_1 \widetilde{U} \end{bmatrix}^{-1}. 
\label{eq:UTS1U}
\end{equation}
On the other hand, by the orthogonality of $U$, 
\begin{equation} (U^T S_1 U)^{-1} = U^T S_1^{-1} U = \begin{bmatrix} \widehat{U}^T S_1^{-1} \widehat{U} & \widehat{U}^T S_1^{-1} \widetilde{U} \\
\widetilde{U}^T S_1^{-1} \widehat{U} \ & \widetilde{U}^T S_1^{-1} \widetilde{U} \end{bmatrix}.
\label{eq:UTSU}
\end{equation}
Thus, $W$ is a product of the leading block of the matrix (to be inverted) on the right-hand side of \eqref{eq:UTS1U}  and the leading block of the matrix on the right-hand side of \eqref{eq:UTSU}. 

For the inversion  of \eqref{eq:UTS1U} we proceed as follows. 
If $M$ is a square nonsingular matrix, $L$ is a square matrix whose dimensions do not exceed those of $M$, and $N_1$ and $N_2$ are  matrices of appropriate dimensions such that the matrix 
$$ R = \begin{bmatrix} M & N_1^T \\ N_2 & -L \end{bmatrix}$$ is well defined, then $R$ is invertible 
if and only if the Schur complement 
 \begin{equation}
S=-(L+N_2 M^{-1} N_1^T) 
\label{eq:S}
\end{equation}
is invertible \cite[Section 3.3]{bgl05}, and its inverse is given by 
\begin{equation} 
R^{-1}= \begin{bmatrix} M^{-1} + M^{-1} N_1^T S^{-1} N_2 M^{-1} & -M^{-1} N_1^T S^{-1}  \\ -S^{-1} N_2 M^{-1}  & S^{-1} 
\end{bmatrix}.
\label{eq:spinv}
\end{equation}

In \cite[Theorem A.1]{Elman99}, the author considered the symmetric case of the Stokes problem with Dirichlet boundary conditions, and applied a symmetric version  of \eqref{eq:spinv} to prove that the minimum eigenvalue of the BFBt-preconditioned Stokes operator (namely, the matrix $\mathcal{M}^{-1} \mathcal{K}$ defined at the beginning of this section, where here $F$ is a vector Laplacian and $G$ is a discrete divergence) is bounded below by 1.

Let us  denote
\begin{equation}
\begin{aligned}
M= & \widehat{U}^T S_1 \widehat{U} \in {\mathbb R}^{p \times p};  &   N_1^T=  \widehat{U}^T S_1 \widetilde{U} \in {\mathbb R}^{p \times (m-p)};  \\
 N_2= &\widetilde{U}^T S_1 \widehat{U} \in {\mathbb R}^{(m-p) \times p}; & L= -\widetilde{U}^T S_1 \widetilde{U} \in {\mathbb R}^{(m-p) \times (m-p)}.
\label{eq:MNL}
\end{aligned}
\end{equation}
Under these assumptions, $S$ is $(m-p)\times (m-p)$ and 
$R$ is $m \times m$.
We then have by \eqref{eq:UTSU} and \eqref{eq:spinv} 
$$ \widehat{U}^T S_1^{-1} \widehat{U} = M^{-1} (I+N_1^T S^{-1} N_2 M^{-1}). $$
From this it readily follows that $W$ defined in \eqref{eq:W} can be written as
 $$ W = I + N_1^T  S^{-1} N_2 M^{-1},$$
where $S$ is as defined in \eqref{eq:S} with  \eqref{eq:MNL} defining $M$, $N_1$, $N_2$, and $L$.

Consider the corresponding eigenvalue problem for $W$
\begin{equation}
(I + N_1^T  S^{-1} N_2 M^{-1} ) x = \gamma x. 
\label{eq:nu}
\end{equation}
If $x=Mz$ where $z \in {\rm ker}(N_2)$, then if this kernel contains nonzero vectors, we have that $\gamma=1$ is an eigenvalue. 

If $m-p<p$, i.e., $p > \frac{m}{2}$, then $N_2$ has more columns than rows and it necessarily has a nontrivial kernel  of dimension $p-(m-p)=2p-m>0$. 
Therefore, in this case $\gamma=1$ in \eqref{eq:nu} is an eigenvalue of multiplicity (at least) $2p-m$.  
 
\subsubsection{Additional observations for the symmetric case of Section~\ref{sec:symD0}}
\label{sec:additional}
We can make additional observations if $S_1$ is symmetric positive definite. In that case,  $U^T S_1 U$ in \eqref{eq:UTS1U} is symmetric positive definite, and so is its leading block, $M=\widehat{U}^T S_1 \widehat{U}$. Furthermore, $N_1=N_2$; let us denote those two matrices by $N$. Therefore, by \cite[Theorem 7.7.6]{horn1990matrix}, the Schur complement
$$ S = \widetilde{U}^T S_1 \widetilde{U} - \widetilde{U}^T S_1 \widehat{U} \left( \widehat{U}^T S_1 \widehat{U}  \right)^{-1} \widehat{U}^T S_1 \widetilde{U}$$
is also symmetric positive definite. Hence, $N_1^T  S^{-1} N_2 M^{-1}$ is similar to the  symmetric positive semidefinite matrix $M^{-1/2} N^T  S^{-1} N M^{-1/2}$, and all 
 eigenvalues of $I + N_1^T  S^{-1} N_2 M^{-1}$ are bounded below by $1$.

The eigenvalues of the corresponding preconditioned matrix $\mathcal{\widehat{M}}_D^{-1} \mathcal{K}$ are, then:
\begin{itemize}
\item $\lambda=1$, of algebraic multiplicity $n-m$.
\item $\lambda=\frac{1 \pm \sqrt{5}}{2}$, (i.e., $\lambda_+\approx 1.618$ and $\lambda_-\approx -0.618$) each of algebraic multiplicity $m-p$.
\item The three distinct eigenvalues 
$2 \cos \left( \frac{\pi } {7} \right) \approx 1.8019, 2 \cos \left( \frac{3 \pi } {7} \right) \approx 0.4450, \\ 2 \cos \left( \frac{5 \pi } {7}\right) \approx -1.2470$,
each of algebraic multiplicity $\max(2p-m,0)$. (Since we assume $m \ge p$, we have $0 \leq \max(2p-m,0) \leq p$.)
\item The remaining eigenvalues are the roots of the cubic polynomials 
$$\lambda^3 - \lambda^2 - (1+\mu_i) \lambda+\mu_i=0, $$
where $ \mu_i$ are the generalized eigenvalues of \eqref{eq:geneig} with (the inverse of) $\widehat{S}_2$  given in \eqref{eq:bfbt} that
 are {\em not equal to 1}; there are $\min(p,m-p)$ such generalized eigenvalues, and since the above polynomials are cubic, there is a total of  $3 \min(p,m-p)$ eigenvalues in this group. 
\end{itemize}

\section{Application: Marker-and-Cell discretization of the Stokes--Darcy equations}
\label{sec:numerical}

	Consider the coupled Stokes--Darcy problem in  $\Omega=\Omega_d \bigcup\Omega_s$,   where  $\Omega_d$ and $\Omega_s$ are two simple rectangular subdomains, connected by $\Gamma$, an interface between two flow regimes.  
	 	A common formulation  is 
	 	\begin{subequations}
	 		\begin{align}
	 			-\nabla\cdot (\kappa \nabla p^d) &=  f^d &  {\rm in}\,\,\Omega_d, \label{eq: Darcy-Laplace-form} \\
	 			-\nu \triangle \bm{u}^s+\nabla p^s&= \bm{f}^s &  {\rm in}\,\,\Omega_s, \label{eq: Dstokes-form1}\\
	 			\nabla\cdot \bm{u}^s&=  0 &   {\rm in}\,\,\Omega_s, \label{eq: Dstokes-form2} 
	 		\end{align}
	 	\end{subequations}
with appropriate boundary conditions.  We refer the reader to \cite{discacciati2009navier} for a detailed description of these equations and their variants.

We solve for $(p^d,\bm{u}^s,p^s)$: $p^d$ is Darcy pressure, $\bm{u}^s=(u^s,v^s)$ is Stokes velocity, and $p^s$ is Stokes pressure. 
	 	The following three interface conditions couple the Stokes and the Darcy equations (respectively, mass conversation, balance of normal forces, and Beavers-Joseph-Saffman condition): 
\begin{equation} 
v ^s= -\kappa \frac{\partial p^d}{\partial y};\qquad 
	 			p^s-p^d = 2 \nu \frac{\partial v^s}{\partial y};  \qquad
	 			u^s  = \frac{\nu }{\alpha} \left(\frac{\partial u^s}{\partial y}+ \frac{\partial v^s}{\partial x}  \right).
\label{eq:interface}
\end{equation}
There are two  physical parameters here: the hydraulic constant $\kappa$ (which appears in this setting in a simplified scalar form), and the viscosity coefficient $\nu.$
	 	We adopt the MAC scheme of \cite{shiue2018convergence}  for discretizing the equations. Uniform staggered grids with meshsize $h=1/n$ are used and the discrete values of $(p^d,u^s,v^s,p^s)$ are placed at different locations: $p^d$ and $p^s$ are evaluated at the cell centers, and the discrete values of $u^s$  and $v^s$ are located at the grid cell faces. Stability and convergence analysis   can be found in \cite{shiue2018convergence}; see also \cite{sun2019stability}.  For details on the discretization of the equations and the three interface conditions \eqref{eq:interface},  
	 	see \cite{shiue2018convergence}.  See  \cite{GH2023, greif2025monolithic} for  iterative solvers for these equations.

We take an example from \cite{luo2017uzawa}, used also in \cite{GH2023}: given $\Omega_s=[0,1]\times [0,1]$ and $\Omega_d=[0,1]\times [-1,0]$, we set
 $u^s =\eta'(y) \cos x, v^s =\eta(y)\sin x,  p^d=0,$ and $ p^s=e^y\sin x,$
 where $\eta(y) = -\kappa -\frac{y}{2\nu}+\left(-\frac{\alpha}{4\nu^2}+\frac{\kappa}{2}\right)y^2.$
The problem is constructed so that this is the analytical solution, and interface conditions are taken into consideration. There are no constraints on $\kappa$ and $\nu$, and we set $\alpha=\nu$.

	 	The resulting system of equations is a double saddle-point system of the form \eqref{eq:system}. Here $A$ is the Darcy operator multiplied by $\kappa$, $D$ is a mildly nonsymmetric and positive definite Stokes operator (the nonsymmetry is a consequence of incorporating the discretized interface conditions), $B$ is a low-rank interface matrix, and $C$ is a discrete divergence operator.

If $n_1$ is the number of meshpoints in one direction, we have in the notation of this paper $n=n_1^2=p$ and $m=2n_1^2-n_1$. 
Thus, the dimensions of the Stokes operator $D$ are larger than those of $A$. 

The rank of the $(2n_1^2-n_1) \times n_1^2 $ 
interface matrix $B$ is $n_1$. Most of the elements of this matrix are zero: its transpose can be written as 
\begin{equation} B^T = \begin{bmatrix} 0 & 0 & 0 \\ 0 & I_{n_1}/h & 0 \end{bmatrix}. 
\label{eq:Bblock}
\end{equation}

\subsection{Validation of the eigenvalue analysis}

\begin{figure}[htb]
		\includegraphics[scale=0.9]{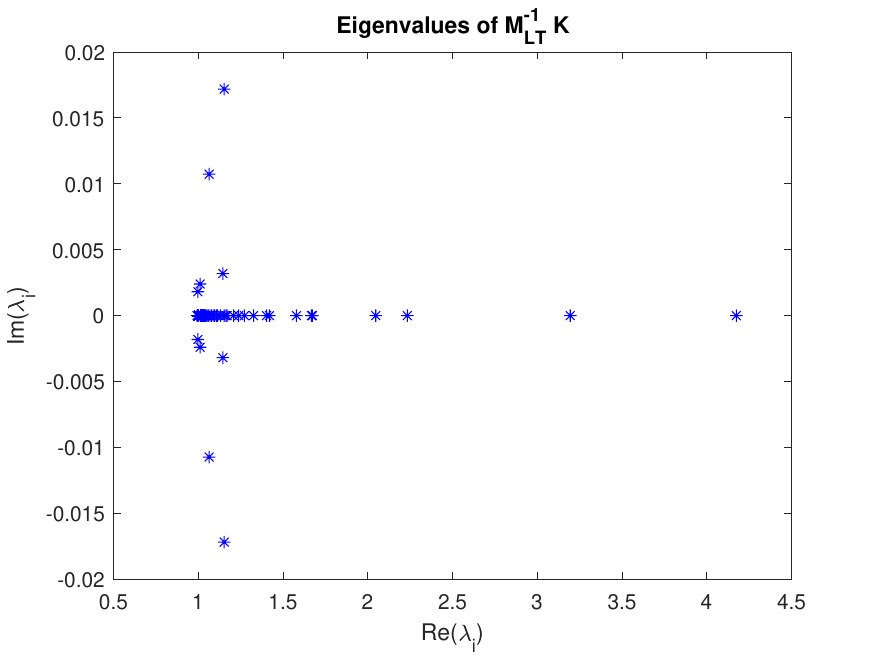}  
		\caption{Stokes-Darcy equations: eigenvalues of $\widehat{\mathcal{M}}_{LT}^{-1} \mathcal{K}$ in the complex plane, with $n=16$, for $\nu=\kappa=1$.  	}
\label{fig:eigT}
\end{figure}

In Figure \ref{fig:eigT} we show  the eigenvalues of the preconditioned matrix $\widehat{\mathcal{M}}_{LT}^{-1} \mathcal{K}$ on a small mesh. We use the exact $A$ and $S_1$, to (partially) validate our eigenvalue analysis. We observe that the eigenvalues are located in the right-half plane on or to the right of the vertical line passing through the eigenvalue $1$. We note that there is a strong concentration of eigenvalues equal to 1 or very close to 1, which is not fully visualized from the plot. 
 Due to the nonsymmetry there are some complex eigenvalues that are hard to characterize. Furthermore, taking practical approximations (see Section \ref{sec:approx}) might generate additional complex eigenvalues, and a complete  eigenvalue analysis is significantly more challenging.
 Nonetheless, the computed eigenvalues are seen to be strongly clustered near 1 and uniformly bounded away from zero, giving rise to fast convergence of the iterative solver.

\begin{figure}[htb]
		\includegraphics[scale=0.9]{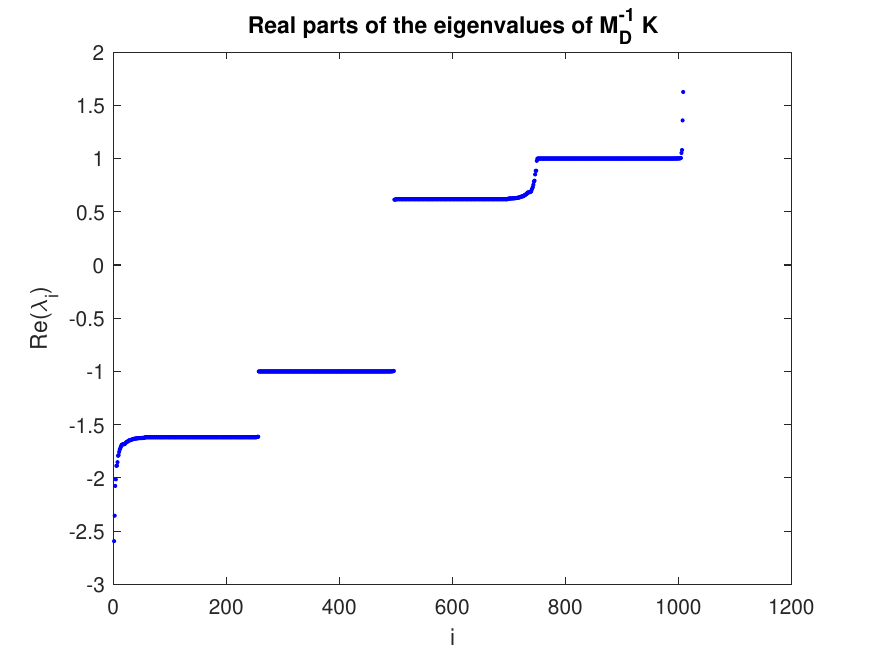}  
		\caption{Stokes-Darcy equations: real parts of the eigenvalues of $\widehat{\mathcal{M}}_D^{-1} \mathcal{K}$, with $n=16$, for $\nu=\kappa=1$.  Some of the eigenvalues are complex but their imaginary parts are smaller than $0.01$ in norm. 	}
\label{fig:eigD}
\end{figure}

In Figure \ref{fig:eigD} we show the real parts of the eigenvalues for the preconditioned matrix $\widehat{\mathcal{M}}_D^{-1} \mathcal{K}.$ The spectral structure closely matches our expectations, as predicted by the analysis in Section \ref{sec:Dne0}. 
While the analysis does not fully hold for this test problem due to the nonsymmetry,  the algebraic multiplicities of the generalized eigenvalues of \eqref{eq:geneig} are equal to the the golden ratio eigenvalues of $\widehat{\mathcal{M}}_D^{-1} \mathcal{K}$ and the algebraic multiplicities of the eigenvalues $1$ and $-1$ match the values predicted by the analysis.

\subsection{Practical approximations of the Schur complements}
\label{sec:approx}

In practice, the Schur complements must be approximated in order for the iterative scheme to converge within a reasonable amount of overall computational work. 
Since the matrix $\mathcal{K}$ is nonsymmetric, we  employ the block-triangular preconditioner $\mathcal{\widehat{M}}_{LT}$ with the BFBt approximation of $S_2$. 
We use an approximation of $S_1$  from \cite{GH2023}:
\begin{equation} 
\widetilde{S}_1 = D+ B A^{-1} B^T \approx D+ \begin{bmatrix} 0 & 0 & 0 \\ 0 & \widetilde{T} & 0 \\ 0 & 0 & 0 \end{bmatrix},
\label{eq:DT}
\end{equation}
where $\widetilde{T}=F^{-T} F^{-1} / h^2$, with $F$ being an incomplete Cholesky with a modest drop tolerance of 0.01 for the Darcy operator $A$. The division by $h^2$ arises from the scale of $B$, as given in \eqref{eq:Bblock}.

We now turn our attention to applying the BFBt approach for $\widehat{S}_2^{-1}$.  We use a diagonal approximation of $A$: $$\widetilde{T}_D=\frac{\tau}{\kappa} I_n,$$ where $\tau=1/3$ is deduced from the discretization; see \cite[Section 4.2.1]{GH2023}.
The matrix $C$ is a discrete divergence operator.
We write it as 
\begin{equation} C = \begin{bmatrix} C_x & C_0 & C_y \end{bmatrix},
\label{eq:Cblock}
\end{equation}
where $C_0=\begin{bmatrix} I_{n_1}/h \\ 0 \end{bmatrix} \in  {\mathbb R}^{n_1^2 \times n_1}.$
Since
$ C S_1 C^T = C D C^T + C B A^{-1} B^T C^T,$
we have
$$ \widehat{S}_2^{-1} = \underbrace{(C C^T)^{-1} C D C^T (C C^T)^{-1}}_{\widehat{S}_{21}} + \underbrace{(C C^T)^{-1} C B A^{-1} B^T C^T (C C^T)^{-1}}_{\widehat{S}_{22}}. $$
For the matrix denoted above by $\widehat{S}_{21}$ we take the  scaled identity $\nu I$ as an approximation, given that $C$ is a discrete divergence operator and $C C^T$ is effectively a scalar Laplacian. 
For $\widehat{S}_{22}$, multiplying $C$ in its block form \eqref{eq:Cblock} by the approximation of $B A^{-1} B^T$  extracted from \eqref{eq:DT} (the matrix $D$ is discarded on the right-hand side) and taking $\widetilde{T}_D$ as the approximation of $\widetilde{T}$, we obtain an approximation
\begin{equation} 
\widetilde{B}_1 = \begin{bmatrix} \widetilde{C}_{11} & 0 & 0 \\ 0 & 0 & 0 \\ 0 & 0 & 0 \end{bmatrix},
\label{eq:B1}
\end{equation}
where $\widetilde{C}_{11} = \frac{\tau}{h^2 \kappa} I_{n_1}.$
We therefore approximate $\widehat{S}_{22}$ by
$\widetilde{S}_{22} = (C C^T)^{-1}  \widetilde{B}_1 (C C^T)^{-1},$
with the above defined choice of $\widetilde{C}_{11}$.
Altogether, we approximate $\widehat{S}_2^{-1}$ by
$$ \widetilde{S}_2^{-1} = \nu I + (C C^T)^{-1}  \widetilde{B}_1 (C C^T)^{-1}. $$

\subsection{Implementation of the preconditioner solves}
\label{sec:implementation}
Since the approximation  appears in  inverse form, it is useful to write the preconditioning operations in terms of inverses, and make sure that no unnecessary inverse operations are done.  The inversion of \eqref{eq:hatMLT} gives
$$ \widehat{\mathcal{M}}_{LT}^{-1} =\begin{bmatrix}
 A^{-1} & & \\
S_1^{-1} B A^{-1} & -S_1^{-1} &\\
-\widehat{S}_2^{-1} C S_1^{-1} B A^{-1} & \widehat{S}_2^{-1} C S_1^{-1} & \widehat{S}_2^{-1}
\end{bmatrix}.$$
In every GMRES iteration, a preconditioner solve requires computing an operator of the form $z=\widehat{\mathcal{M}}_{LT}^{-1} r$, for a given vector $r$.
Suppose
$$   r = \begin{bmatrix} r_1 \\ r_2 \\r_3 \end{bmatrix}.$$
In practice, we replace $S_1^{-1} $ by $\widetilde{S}_1^{-1}$ and $\widehat{S}_2^{-1}$ by $\widetilde{S}_2^{-1}$, as described in Section~\ref{sec:approx}, and an effective approximation of $A^{-1}$, denoted by $\widetilde{A}^{-1}$, is also needed.
Then
\begin{eqnarray*}
z  & = &
\begin{bmatrix}
 \widetilde{A}^{-1} & & \\
\widetilde{S}_1^{-1} B \widetilde{A}^{-1} & -\widetilde{S}_1^{-1} &\\
-\widetilde{S}_2^{-1} C \widetilde{S}_1^{-1} B \widetilde{A}^{-1} & \widetilde{S}_2^{-1} C \widetilde{S}_1^{-1} & \widetilde{S}_2^{-1} \end{bmatrix}
\begin{bmatrix} r_1 \\ r_2 \\r_3 \end{bmatrix} \\ 
& = &
\begin{bmatrix}
\widetilde{A}^{-1} r_1 \\
\widetilde{S}_1^{-1} (B \widetilde{A}^{-1} r_1 -r_2)\\
-\widetilde{S}_2^{-1} (C \widetilde{S}_1^{-1} (B \widetilde{A}^{-1} r_1 -  r_2) - r_3)
\end{bmatrix}.
\end{eqnarray*}
Therefore,  $t_0=\widetilde{A}^{-1} r_1, \ t_1=B t_0$, and $t_2=\widetilde{S}_1^{-1} t_1$ are computed  once and stored, and inversions or multiplications by inverted factors of approximations of $A$, $S_1$, and $S_2$ are also computed only once. 

In a purely iterative environment, some of the above quantities are typically computed approximately using iterative inner solvers.  It is also necessary to consider inner iterations for solving linear systems associated with $C C^T$. Given the close relationship of this matrix to a scalar Laplacian, this can be done by a fast multigrid solver.

\subsection{Performance of the iterative solver}

In our {\sc Matlab} implementation, run on an i7 3.4GHz Intel processor with 128GB of RAM, we have used the approximations described  in Section~\ref{sec:approx} and the implementation of preconditioner solves described in Section~\ref{sec:implementation}. However, we did not implement fast solvers or used approximations for inversions of $A$ and $C C^T$. 

Tables \ref{tab:nu1e0inexact}--\ref{tab:nu00001e0inexact} present results for a few choices of $\nu$ and $\kappa$. The dimensions of the linear systems are quadrupled each time the mesh is refined by a factor of 2. For $n_1=32$ the system dimensions are $4,064 \times 4,064$, and for the largest value of $n_1$ in the experiments, $n_1=512$, the dimensions are $1,048,064 \times 1,048,064$. Ideally, we hope to see minimal sensitivity of the iteration counts to the mesh size and the values of the physical parameters. In practice, we see that there is a loss of speed of convergence as the mesh is refined, but it is rather slow, and similarly, convergence appears to be relatively robust with respect to the values of the physical parameters. For $n_1=512$ and $\nu=10^{-4}$  we have observed stagnation; see Table \ref{tab:nu00001e0inexact}. We note that even in these cases the iterations  trended in the direction of a decrease of the residual norm.  
For small values of $\kappa$, there is a stronger change in the iteration counts as the mesh is refined. 

Generally speaking,   increased numerical difficulties are expected for small values of the physical parameters. 
In particular, the loss of scalability may  manifest itself in an increased computational cost of inner iterations for the components of the block preconditioner. 

The CPU times appear to reasonably scale with respect to the linear system sizes, although they do not exhibit optimal scalability as the mesh is refined; we speculate that this is primarily due to the different optimization levels of builtin {\sc Matlab} functions and the non-optimal implementation of the inner solves.

\begin{table}
 \caption{Iteration counts and CPU times (in seconds) for GMRES(20) and  the preconditioner $\widehat{\mathcal{M}}_{LT}$ (using the  approximations described in Section~\ref{sec:approx}) with $\nu=1$ and varying $n_1$ and $\kappa$.}
\centering
\begin{tabular}{c|ccccccccc}
\hline
\diagbox{$n_1$}{$\kappa$} &   $10^0$  & $10^{-2}$  & $10^{-4}$ &  $10^{-6}$    \\ \hline     
32 &     19 (0.022)  &  17 (0.021)  &  15 (0.019) &  12 (0.018)  \\
64 &     20 (0.12)   &  18 (0.11)  &  17 (0.10) &  14 (0.092)  \\
128 &   21 (0.89) &  19 (0.78)  &  19 (0.79)  &  16 (0.67) \\
256 &   22 (5.05) &  20 (4.55)  &  23 (5.24)   &  19  (4.33) \\
512 &   23 (29.71) &   20 (25.06)  &   27 (33.03)  &   25 (31.61) \\
 \hline
\end{tabular}\label{tab:nu1e0inexact}
\end{table}

\begin{table}
 \caption{Iteration counts and CPU times (in seconds) for GMRES(20) and  the preconditioner $\widehat{\mathcal{M}}_{LT}$ (using the  approximations described in Section~\ref{sec:approx}) with $\nu=0.01$ and varying $n_1$ and $\kappa$.}
\centering
\begin{tabular}{c|ccccccccc}
\hline
\diagbox{$n_1$}{$\kappa$}   &  $10^0$  & $10^{-2}$  & $10^{-4}$ &  $10^{-6}$    \\ \hline     
32 &     16 (0.021)   &  12 (0.015)   &  12  (0.015)  &  13 (0.016)  \\
64 &     17  (0.11) &  14 (0.094)  &  12 (0.082)  &  14  (0.091) \\
128 &   17 (0.76)  &  15  (0.67)  &  13 (0.60)  &  16 (0.72)  \\
256 &   17 (4.11) &  17 (4.32) &  15 (4.02)  &  17 (4.31) \\
512 &   17 (23.65) &   19 (26.45)  &   17 (25.10) &   17 (24.30) \\
 \hline
\end{tabular}\label{tab:nu001e0inexact}
\end{table}  

\begin{table}
 \caption{Iteration counts and CPU times (in seconds) for GMRES(20) and the preconditioner $\widehat{\mathcal{M}}_{LT}$  (using the approximations described in Section~\ref{sec:approx})  with $\nu=0.0001$ and varying $n_1$ and $\kappa$. In the last row, `S' stands for stagnation: the iterations seemed to converge but there were stagnation error messages indicating that two consecutive iterates were the same.}
\centering
\begin{tabular}{c|ccccccccc}
\hline
\diagbox{$n_1 $}{$\kappa$}   &  $10^0$  & $10^{-2}$  & $10^{-4}$ &  $10^{-6}$    \\ \hline     
32 &     12 (0.017)  &  10 (0.014) &  12 (0.015) &  13 (0.016)  \\
64 &     14 (0.094) &  10  (0.074) &  13 (0.093)  &  15 (0.097)  \\
128 &   15  (0.73) &  10 (0.53)  &  13 (0.68) &  17  (0.78) \\
256 &   17  (4.52) &  12 (3.74)   &  13 (3.75)  &  19 (5.10)  \\
512 &   19 (28.82) &   S &   S &   S \\
 \hline
\end{tabular}\label{tab:nu00001e0inexact}
\end{table}

\section{Concluding remarks} 
\label{sec:conc}

 We have shown the important role that the spectrum of $\widehat{S}_2^{-1} S_2$ plays in addressing the question to what extent the quality of the approximation of $S_2$ affects the convergence  of the iterative solver.

The block-triangular BFBt/least-squares commutator has been successfully applied to the Navier-Stokes equations, which feature block $2 \times 2$ nonsymmetric systems. In this work we have adapted this approach to the approximation of $S_2$ in a double saddle-point system. A key to the effectiveness of this approach is that $C C^T$ should be easy to invert; this is indeed the case in the Stokes-Darcy equations that we have used as an example. 

Extending the use of this preconditioner, making it as practical as possible in the face of the challenge in approximating $S_2$, and obtaining additional analytical observations on the complex eigenvalues of nonsymmetric preconditioned double saddle-point systems, are desirable goals for future work.

\section*{Acknowledgments} I am deeply grateful to two  referees whose  helpful comments and suggestions have greatly improved  the quality of this manuscript.

\bibliographystyle{plain}
\bibliography{g25-ref}

\end{document}